\documentclass[11pt,twoside, final]{amsart}  % AMS article style
 \usepackage{a4wide}
  \copyrightinfo{0}{Lorestan University}
 \pagespan{1}{\pageref*{LastPage}}
\usepackage{etoolbox,lastpage}
 \usepackage{booktabs} % for creating desired tables
  \usepackage{amsmath,amsthm,amscd,amsfonts,amssymb}
   \usepackage{graphicx}	  % For eps figures	
    \usepackage{color}        % For color text: \color command
     \usepackage[colorlinks]{hyperref}
      \usepackage{fancyvrb,enumerate,array}
       \usepackage{fvextra,wrapfig,lipsum}
        \usepackage{multirow,multicol}
        \usepackage{comment}
        \usepackage{bigints}
         \commby{}

    \theoremstyle{definition}

       \newtheorem{Q}{Question}

\newtheorem*{Q*}{Question}
\newtheorem*{Thm*}{Main Theorem}
\newtheorem*{THEO*}{Theorem[Khanehdani-Swua]\cite{Khanehdani}}
\newtheorem*{Proof}{Proof of the Main Theorem}
\theoremstyle{remark}
 \numberwithin{equation}{section}
\begin{document}

	\begin{comment}
	
\begin{table}[t]
	   \vspace{-1.5cm}
	   \noindent
	   \hspace{-.5cm}	
	   \begin{tabular}{llll}
			\includegraphics[width=1.5cm]{logo.eps} &
			\begin{tabular}{ll}\vspace{-1.6cm}\\
			   \hspace{-2cm}{\Large \bf Mathematical Analysis} & \\[2ex]
			   \hspace{-2cm}{\Large \bf	\& Convex Optimization } &
			\end{tabular}
			
			 & \qquad \qquad   \qquad   &
			
			\begin{tabular}{l}\vspace{-1.5cm}\\
		  	   {\large Vol. xxxxxx }\\
			   {\large https:$\backslash \backslash$ maco.lu.ac.ir}\\
	    	   {\large  DOI: 0000000000000}		
			\end{tabular}\\		
      \hline  \\[-2ex]
	  {\large \bf Research Paper}
      \end{tabular}
\end{table}
\end{comment}

%% The title of the paper goes here.  Edit your title.
\title[Complex Limit Cycle]{ A Complex Limit Cycle not Intersecting the Real Plane  }
%% Now edit the following to give First Author name and address:
%% $^*$ for the corresponding author.

\author[A. Taghavi]{Ali Taghavi}
\address[Ali Taghavi]{Department of Mathematics,
	Qom University of Technology,Qom, Iran.}
\email{taghavi@qut.ac.ir}

%\author[f2. last-name2]{first-name2 last-name2$^*$}
%\address[first-name2 last-name2]{Complete Address}
%\email{SecondEmailAdress@some.ir}
%% If there are three of more authors they are added in the obvious way.
 %\thanks{We thank Lo\"{\i}c Teyssier for his very helpful comments and suggestions}
%------------------------------------------------------------------------------------%
%%
%% Use the following command to make the title for the paper.
%
 %\CoverPage
%\date{Received:  , Accepted:  .}
\maketitle
%
%%% The following environment is needed for the abstract.
%%%

\begin{abstract}
We give a precise example of a polynomial vector field on $\mathbb{R}^2$ whose corresponding singular holomorphic foliation of $\mathbb{C}^2$ possesses a complex limit cycle which does not intersect the real plane $\mathbb{R}^2$. \\ \\


\end{abstract}

\section{\bf Introduction}
Every  polynomial vector field $\begin{cases} \frac{dz}{dt}=P(z,w)\\ \frac{dw}{dt}=Q(z,w) \end{cases}$   defines  a  singular foliation by holomorphic curves on $\mathbb{C}^2$.   The qualitative study of the behaviour of these curves plays a crucial role  in the area of holomorphic dynamical system.  Every nonsingular curve of the complex equation is considered as  a leaf of the corresponding singular foliation. In the real dynamical system a limit cycle is a closed orbit whose  Poincare return  map is not  equal to the identity map. For introduction on limit cycle phenomena see \cite[page 250]{hirsch-smale}. The number of limit cycles of polynomial vector fields  is the main object of the second part of the Hilbert 16th Problem. For details on Hilbert 16th problem see \cite{Ilya}.  However in the complex setting a complex limit cycle is a complex leaf  which contains a closed curve $\gamma$ whose holonomy map  is not equal to the identity map. The concept of  holonomy is an immediate complex generalization of the  Poincare return map. Historically, the idea of consideration
 of such kind of foliations arising from algebraic vector fields goes back to the innovative idea by I. Petrovskii and E. Landis who suggested  considering a real polynomial vector field as a complex vector field. So according to this idea we consider a    real limit cycle as  the intersection of a complex limit cycle with the real plane. For the seminal paper of I. Petrovskii and E. Landis see \cite{PL}. For a fundamental and deep contribution to this innovative idea of Petrovski-Landis, see the influential paper by Yu. Ilyashenko in \cite{Yu}.

In this paper we consider the following question:

\begin{Q*}
 Is there a polynomial vector field on $\mathbb{R}^2$ whose corresponding singular holomorphic foliation on $\mathbb{C}^2$  possess a complex limit cycle $L$   which does not intersect the real plane $\mathbb{R}^2$?
 \end{Q*}

By real plane we mean $\mathbb{R}^2\subset \mathbb{C}^2$. That is  $$\{(z,w)\in \mathbb{C}^2\mid \Im(z)=\Im(w)=0\}.$$

A motivation for this question is the following:  Every real limit cycle is the intersection of a complex limit cycle with the real plane. But in the above question we somehow search for a precise example of an  invisible limit cycle, namely a complex limit cycle which does not intersect the real plane.

We shall prove that the answer to the above question is affirmative as it is indicated in the following theorem. In this theorem by $z',w'$ we mean $dz/dt, dw/dt$ respectively where $t\in\mathbb{C}$ is a complex time-parameter.

\begin{Thm*}
The polynomial vector field \begin{equation}\label{SS} \begin{cases} z'=w+z(z^2+w^2+1)\\ w'=-z+w(z^2+w^2+1) \end{cases} \end{equation} has $z^2+w^2+1=0$ as a complex limit cycle which does not intersect the real plane $\mathbb{R}^2$.

\end{Thm*}

\section*{Preliminaries and Proofs} A $k$ dimensional  foliation of a complex  manifold $M$ is a decomposition of $M=\bigcup L_{\alpha}$ where each $L_{\alpha}$ is an immersed k dimensional holomorphic submanifold of $M$. Moreover for every $p\in M$ there is a neighborhood  $U$ and a bioholomorphism $\phi:U\to V\times W\subset \mathbb{C}^k\times \mathbb{C}^{n-k}$,  where   $V$  and $W$ are open poly discs   in $\mathbb{C}^k$ and $ \mathbb{C}^{n-k}$ respectively. Moreover the bioholomorphism $\phi$ maps each connected component of $L_{\alpha} \cap U$ to a plaque $V\times\{c\}$ for a fixed $c\in W$. The triple $(U,V\times W,\phi)$ is called  a foliation chart.  An important example of a holomorphic foliation of $\mathbb{C}^2$ arises from a non vanishing holomorphic vector field. More precisely, to a vector field \begin{equation}P(z,w)\partial/\partial z+Q(z,w)\partial /\partial w \end{equation} on $\mathbb{C}^2$ one associate a differential 1-form \begin{equation}\omega=-Q(z,w)dz+P(z,w)dw .\end{equation} Then   leaves of the foliation defined by $\omega=0$ are the integral curve of the initial vector field. For more details on real and holomorphic foliations see \cite{Lawson} and \cite{Zakeri}. \\

  Assume that $L$ is a leaf of a holomorphic foliation. Let $\gamma:[0,1]\to L$ be a continuous curve in $L$. At points $\gamma(0)$ and $\gamma (1)$ we choose two $n-k$ dimensional transversal sections $S_0$ and $S_1$, respectively. They are $n-k$ dimensional complex  local sections which are transversal to the leaf $L$. These transversal sections play the role of Poincare sections in the real dynamical system. Then holonomy map along $\gamma$ is a holomorphic function $h:S_0\to S_1$  with the property that $z\in S_0$ and $h(z)\in S_1$ lie on the same leaf of the foliation. This concept is the complex analogy of the Poincare map in the real context. In the real case we  have a regular orbit $\gamma$ with two transversal sections $S_0, S_1$. Then the Poincare map is a well defined map from $S_0$ to $S_1$ which is naturally defined via flow of the vector field. In the Holomorphic case we simulate the same situation but via foliation charts moving along $\gamma$  since the flow has no a well behaved  definition in the holomorphic case. For  details on definition and geometric description of the  holonomy map of a foliation  see \cite{Zakeri}.
For a leaf $L$ of foliation we fix a base point $p$. For every loop $\gamma\subset L$ based at $p$ one may consider the holonomy map $h_{\gamma}$. This gives an obvious group structure on the set $\{h_{\gamma}\mid \gamma\in  \pi_1(L, p)\} $, where $\pi_1(L,p)$ is the fundamental group of the leaf $L$ associated to loops based at point $p$. This group is called the monodromy group of $L$. We denote it by $G(L)$. So there is a natural group homomorphism from $\pi_1(L)$ into $G(L)$.  A  complex limit cycle is  a  leaf $L$ whose monodromy  group $G(L)$ is a non trivial group.\\

Every differential 1-form $\omega=-Q(z,w)dz+P(z,w)dw$ defines a 1 dimensional singular holomorphic foliation of $\mathbb{C}^2$. The leaves of this foliations are the integral curves of \begin{equation}\label{EQ2}   \begin{cases} z'=P(z,w)\\w'=Q(z,w) \end{cases} \end{equation} A natural question is that in this special case that the foliation is given by a global differential 1-form, how can we compute the holonomy of a typical leaf $L$  of the foliation along a closed  curve $\gamma$ lying in $L$?  The following beautiful theorem by Khanehdani-Swua gives us a formula for the first variation of the holonomy map. Before we state the theorem Khanehdani-Swua  we give a very short introduction based on \cite[Khanehdani-Swua]{Khanehdani}:\\
 Suppose that  the foliation is defined via equation $$\omega=-Q(z,w)dz+P(z,w)dw=0.$$ Assume that $L$ is a leaf of the foliation and $\gamma$ is an arbitrary closed curve on $L$. Then we can write $$d\omega=\alpha\wedge \omega,$$  where $\alpha$ is a multi-valued 1-form in a neighbourhood of $\gamma$  and the  restriction of $\alpha$ to each leaf is   single-valued. Assume that $h_{\gamma}$ is the holonomy of $L$ associated to $\gamma$ defined in a neigbourhood of origin in $\mathbb{C}$ which models a local section transversal to the leaf. Then we have:

\begin{THEO*}
The first variation of the holonomy is given by  $h_{\gamma}'(0)=\exp{\bigintsss_{\gamma} \alpha}$.

\end{THEO*}

The proof of our main theorem is based on usage of the above theorem:

\begin{Proof}
First we show that the algebraic curve $L: z^2+w^2+1=0$ is  an integral curve of the equation \ref{SS}.  We compute  the derivation $z^2+w^2+1$ along the vector field \ref{SS}: \begin{equation*} \begin{aligned}(z^2+w^2+1)'=2zz'+2ww'=\\&2z(w+z(z^2+w^2+1))+2w(-z+w(z^2+w^2+1))=\\ &(2z^2+2w^2)(z^2+w^2+1).\end{aligned}\end{equation*} This obviously shows that the algebraic curve $z^2+w^2+1=0$ is an integral curve of equation \ref{SS}. It is a non singular solution which does not intersect  the real plane $\mathbb{R}^2\subset \mathbb{C}^2$. Now we prove that $L$ is a complex limit cycle. To prove this we give a closed curve $\gamma$ on $L$ such that the holonomy map $h_{\gamma}$ is different from the identity map. For a differential 1-form $\omega=-Qdz+Pdw$ we have $d \omega=\alpha \wedge \omega$, where \begin{equation}\label{aaa} \alpha=(\frac{P_z+Q_w}{P^2+Q^2})(Pdz+Qdw).\end{equation} We apply this formula to \begin{equation} \begin{cases}P(z,w)=w+z(z^2+w^2+1)\\ Q(z,w)=-z+w(z^2+w^2+1)\end{cases}  \end{equation} We shall  introduce a closed curve $\gamma\subset L$ for which $\int_{\gamma} \alpha \neq 2k\pi i,\quad \forall k\in \mathbb{Z}$. This would show that the holonomy map $h_{\gamma}$ is different from the identity map since $h_{\gamma}'(0)\neq 1$. Notice that the restriction of $\alpha$ in \ref{aaa} to $L$ is equal to $2(wdz-zdw)$.
\newpage
Because $P_z+Q_w=2(z^2+w^2+1)+2z^2+2w^2 \equiv  -2$ on $L$. Similarly  $P^2+Q^2\equiv -1$ on $L$. So $\alpha$ restricts to $wdz-zdw$ on curve $L$. To compute  $\bigintss_{\gamma} \alpha$ we give a global parametrization of $L: z^2+w^2+1=0$  then we apply the change of variable formula in integration:
\begin{equation}\phi:\mathbb{C}\setminus\{0\} \to L\quad \phi(t)=(z(t), w(t)),\end{equation} where
$$z(t)=\frac{i}{2} (t+\frac{1}{t}),  w(t)= \frac{1}{2}(t-\frac{1}{t})).$$
Note that the global parametrization $\phi$ of $L$ gives us an embedding of $\mathbb{C}\setminus\{0\}$ into $\mathbb{C}^2$ whose image is whole $L$.
Moreover we have \begin{equation*} \begin{aligned}\phi^*(wdz-zdw)=w(t)d(z(t))-z(t)d(w(t))=\\&\frac{i}{4} \left((t-\frac{1}{t})(1-\frac{1}{t^2})-(t+\frac{1}{t})(1+\frac{1}{t^2})\right).  \end{aligned} \end{equation*}

The closed curve $\gamma$ in $L$ whose holonomy map $h_{\gamma}$ is nontrivial is the curve $$\phi(e^{i\theta}),\quad \theta \in [0,2\pi].$$ Put $\beta(\theta)=e^{i\theta}, \theta \in [0,2\pi]$, the unit circle in $\mathbb{C}\setminus \{0\}$. So we have $\phi(\beta)=\gamma$. For this closed curve $\gamma$ we have $\int_{\gamma} \alpha=\int_{\beta} \phi^* (\alpha).$ This is a consequence of change of variables in integration. As we computed in the above lines, the restriction of $\alpha$ to $L$ is $\alpha=2(wdz-zdw)$  so we can write:

\begin{equation*} \begin{aligned}\bigintssss_{\gamma} \alpha=\bigintssss_{\beta} \phi^* (\alpha)=\\&\frac{i}{4}\bigintsss_{\beta}\left((t-\frac{1}{t})(1-\frac{1}{t^2})-(t+\frac{1}{t})(1+\frac{1}{t^2})\right) dt=\\&\frac{i}{4}\bigintssss_{\beta} (\frac{4}{t}+\frac{2}{t^3})dt=-2\pi, \end{aligned} \end{equation*} by residue formula in complex function theory, see \cite{Ahlfors}

 So according to theorem by Khanehdani-Swua we see that  the first variation of the holonomy is $e^{-2\pi}$ which is different from $1$.  So the holomomy map is not the identity map. Thus $L:z^2+w^2+1=0$ is  a complex limit cycle  which obviously does not intersect the real plane $\mathbb{R}^2$.  This completes the proof of the main theorem.

\end{Proof}

\section*{Discussions and Further researches}
In this paper, inspired by the relation  between real and complex limit cycle theory ,  we introduced   a precise example of a polynomial vector field on $\mathbb{R}^2$  whose corresponding complex foliation  possesses a complex limit cycle not intersecting the real plane. In fact  we observed that the algebraic curve $L:z^2+w^2+1=0$ is a complex\newline  limit cycle for the real coefficient algebraic vector field
 $$ \begin{cases} z'=w+z(z^2+w^2+1)\\ w'=-z+w(z^2+w^2+1) \end{cases} $$
     The algebraic property of $L$ in our example is a motivation to ask the following question:

 \begin{Q}
 Is there an algebraic vector field on $\mathbb{R}^2$ whose corresponding singular foliation of $\mathbb{C}^2$ possess a complex limit cycle $L$ which does not intersect the real plane $\mathbb{R}^2$ and $L$ is not an algebraic curve?
 \end{Q}

 Moreover the construction of a complex limit cycle not intersecting the real plane naturally arises the following question:

 \begin{Q}
 Is there a polynomial vector field $X$ on $\mathbb{R}^2$ such that every complex limit cycle of the corresponding foliation of $\mathbb{C}^2$ must necessarily intersect the real plane $\mathbb{R}^2$? In this question we do note require this intersection would be in form of a closed curve but we merely require the intersection would be a non empty set.
 \end{Q}

 \section*{ACKNOWLEDGMENT}

I would like to thank Lo\"{\i}c Teyssier for his very helpful  suggestions. I would like also to thank the referees for their very helpful comments  on the manuscript. This work was supported by the Qom University of Technology.

\end{document}